\documentclass[12pt]{article}
\usepackage{bm}

\usepackage{amssymb,amsmath,amsthm,latexsym}

\usepackage{epsfig}

\def\be{\begin{equation}}
\def\ee{\end{equation}}
\def\bea{\begin{eqnarray}}
\def\eea{\end{eqnarray}}
\def\pt{\partial}
\def\ffi{\varphi}

\def\al{\alpha}
\def\gm{\gamma}
\def\Gm{\Gamma}
\def\Th{\Theta}
\def\eps{\varepsilon}
\def\Dt{\Delta}

\def\La{\Lambda}

\def\const{\mbox{const}}

\def\mod{\mbox{mod}}
\def\ln{\mbox{ln}}

\def\dd{\mbox{d}}

\def\calP{{\cal P}}
\def\calH{{\cal H}}
\def\calS{{\cal S}}

\begin{document}

\title{Stability islands in domains of separatrix crossings in
 slow-fast Hamiltonian systems
}

\vskip 15mm

\author{Anatoly Neishtadt$^{1}$, Carles Sim\'o$^{2}$, Dmitri Treschev$^{3}$,
and Alexei Vasiliev$^{1}$
\\
\\
$^{1}$ {\small Space Research Institute, Profsoyuznaya 84/32, Moscow 117997, Russia,}
\\
$^{2}$ {\small Departament de Matem\`atica Aplicada i An\`alisi, Univ. de Barcelona,}\\
{\small Gran Via, 585, Barcelona 08007, Spain, }\\
$^{3}$ {\small V.A. Steklov Mathematical Institute, Gubkina 8, Moscow 119991, Russia.} }

\date{}
\maketitle

\begin{abstract}
We consider a 2 d.o.f. Hamiltonian system with one degree of
freedom corresponding to fast motion and the other corresponding
to slow motion. The ratio of the time derivatives of slow and fast
variables is of order $0<\eps \ll 1$.  At frozen values of the
slow variables there is a separatrix on the phase plane of the
fast variables and there is a region in the phase space (the
domain of separatrix crossings) where the projections of phase points
onto the plane of the fast variables repeatedly cross the
separatrix in the process of evolution of the slow variables.
Under a certain symmetry condition we prove existence of many (of
order $1/\eps$) stable periodic trajectories in the domain of the
separatrix crossings. Each of these trajectories is surrounded by
a stability island whose measure is estimated from below by a
value of order $\eps$. Thus, the total measure of the stability
islands is estimated from below by a value independent of $\eps$.
The proof is based on an analysis of the asymptotic formulas for the
corresponding Poincar\'e map.

\end{abstract}

%\ams{34E10, 37J40}

%\maketitle

\section{ Introduction}\label{intro}

Many problems in the theory of charged particles' motion, the theory of propagation of short-wave excitations, and in
celestial mechanics can be reduced to the analysis of 2 d.o.f. Hamiltonian systems with fast and slow variables
\cite{BZ,GZ,Wis}. One degree of freedom corresponds to fast variables, and the other corresponds to slow
variables. The ratio of the time derivatives of slow and fast variables is of order $0<\eps \ll 1$. To describe the
dynamics in such systems one can use the adiabatic approximation constructed as follows.

Consider the fast system, i.e. the system for the fast variables
at frozen values of the slow variables. This is a 1 d.o.f.
Hamiltonian system involving the slow variables as parameters.
Assume that for a range of values of the slow variables there is a
region filled with closed trajectories on the phase portrait of
the fast system. Then one can introduce ``action-angle'' variables
in the fast system \cite{Arn}. The ``action'' variable of the fast
system is an adiabatic invariant (i.e. an approximate first
integral) of the complete system: its value oscillates with
amplitude $\sim \eps$ on time periods of order $\sim 1/\eps$ along
a phase trajectory. To describe approximately the variation of the
slow variables, one should average the rates of their variation
over the ``angle'' variable of the fast motion. This approximation
of the real motion is called the adiabatic approximation
\cite{AKN}. The obtained 1 d.o.f. Hamiltonian system for the slow
variables, involving the ``action'' variable of fast motion as a
parameter, is called the slow system.

Assume that for some range of values of the slow variables the ``action'' variable phase
trajectories of the slow system are closed. Then in this range the motion can be
described with V.I.Arnold's theorem on perpetual adiabatic invariance \cite{Arn63}: the
``action'' value along a trajectory perpetually undergoes only small, of order $\eps$,
oscillations (under a certain generality condition). The considered region of the phase
space is filled, up to a small residual set, with invariant tori, which are
$O(\eps)$-close to the invariant tori of the system of the adiabatic approximation.

We shall consider cases in which there are separatrices on the fast system's phase portrait (see Figure
\ref{fastplane}), and there is a region in the phase space (the domain of separatrix crossings) where the
projections of phase points onto the plane of the fast variables repeatedly cross the separatrix in the
process of evolution of the slow variables. In this case the above described construction of the adiabatic
approximation needs some modification, and the above formulated assertions on the adiabatic
approximation accuracy are not applicable.

 We shall assume that a certain symmetry condition is valid: the
areas inside the separatrix loops in Figure \ref{fastplane} are equal.

The main result of the present paper is as follows. Under certain generality conditions in the domain of
separatrix crossings on every energy level there exist many, of order $1/\eps$, stable periodic trajectories
of period $\sim 1/\eps$. Each of these trajectories is surrounded by a stability island, and the measure of
this island is estimated from below by a value of order $\eps$. Thus, the total measure of the stability
islands is estimated from below by a quantity which is independent of $\eps$. A stability island is a domain
on an energy level bounded by a two-dimensional invariant torus. A stability island contains a discrete
family of invariant tori contractible to the periodic trajectory. Introduce a ``modified action'' equal to the
``action'' for points inside the separatrix loops in Figure \ref{fastplane}, and equal to one half of the
``action'' for the other points. This ``modified action'' is a perpetual adiabatic invariant of the motion inside
a stability island: its value along a phase trajectory perpetually undergoes only oscillations with amplitude
of order $\eps$. Therefore, the stability islands are also islands of perpetual adiabatic invariance.

We note that the existence of stability islands with total measure that is not small with
$\eps$ in the domain of separatrix crossings is quite unexpected. Visually, in many
problems, this domain looks like a region of dynamical chaos (see \cite{CB,NS04}). It was
shown in \cite{EE} that a single stability island cannot have a measure larger than
$O(\eps)$.

The existence of stability islands with total measure not small with $\eps$ in the domain of separatrix crossings
was established in \cite{NST97,NST2001} in the case of a Hamiltonian system with one degree of freedom and
the Hamiltonian function slowly periodically depending on time. Here we generalize this result to 2 d.o.f.
systems. Like in \cite{NST97,NST2001}, the proofs are based on the study of asymptotic formulas for the
corresponding Poincar\'e map. In \cite{NST97,NST2001} these formulas were constructed with the use of
asymptotic expressions for the jump of the adiabatic invariant at a separatrix crossing in systems with one
and a half d.o.f. \cite{T,CET,N86} and for the variation of the ``angle'' variable between separatrix crossings
in such systems \cite{CS}. Below we use analogous formulas for systems with two d.o.f. \cite{N87,NV05}. To
find linearly stable periodic trajectories we look for linearly stable fixed points of the Poincar\'e map.
The results on Lyapunov stability of the periodic trajectories (the fixed points) and on the existence of
invariant tori surrounding the periodic trajectories (invariant curves around the fixed points) are provided
by the Kolmogorov-Arnold-Moser (KAM) theory. An example of application of the present theory to the motion of a
charged particle in the parabolic model of magnetic field in the Earth magnetotail is given in \cite{NSTV}.

\section{ Formulation of the results}\label{results}

Consider a two degrees of freedom Hamiltonian system with Hamiltonian $H = H(p,q,y, x)$, where $q, \eps^{-1}
x$ are coordinates, and $p, y$ are canonically conjugated momenta, $\eps >0$ is a small parameter, $H \in
C^{\infty}$. The corresponding equations of motion are:

\be
\dot p = -\frac{\pt H}{\pt q}, \;\; \dot q = \frac{\pt H}{\pt p}, \;\; \dot{y} = -\eps \frac {\pt
H}{\pt x}, \;\; \dot{x} = \eps \frac{\pt H}{\pt y}.
\label{2.0}
\ee
The variables $p,q$ are called fast variables, while $y, x$ are called slow variables. The Hamiltonian system
for $p,q$ at $(y, x)=\const$ is called fast (or unperturbed) system.

\begin{figure}[ht]
\hspace*{44mm}\epsfig{file=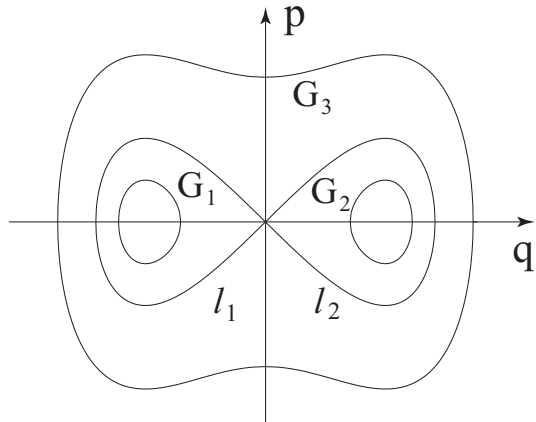,width=70mm} \vspace*{-5mm} \caption{Phase portrait of the fast
system.}
\label{fastplane}
\end{figure}

Assume that at all considered values of the slow variables the phase portrait of the fast system looks as in
figure \ref{fastplane}. On this phase portrait, there exist a non-degenerate saddle point $C$ and
separatrices $l_1, l_2$ . These separatrices divide the $(p,q)$-plane into domains $G_i = G_i(y, x) , \,
i=1,2,3$. Denote the value of $H$ at the saddle point as $h_s = h_s(y, x)$, and introduce $E = E(p,q,y, x) =
H - h_s$. On the separatrices $E=0$. We suppose that $E>0$ in $G_3$ and $E<0$ in $G_1$ and $G_2$. The areas of
the domains $G_1, G_2$ are $S_1, S_2$ respectively, and $S_3 = S_1 + S_2$; $S_i = S_i(y, x)$ . We will
make the following symmetry assumption:

A) The areas of $G_1, G_2$ are equal: $S_1(y, x)=S_2(y, x)=S(y, x)$.

Put $\Th =\Th(y, x)= \{S,h_s\}$, where $\{\cdot,\cdot\}$ is the Poisson bracket with respect to the variables
$(y, x): \; \{f,g\}=f'_{x} g'_{y} - f'_{y} g'_{x}$.

In the unperturbed system, one can introduce canonical ``action-angle'' variables $(I,\varphi)$ separately in
each $G_i$. The ``action'' variable $I=I(p,q,y, x)$ in the unperturbed system involves $y, x$ as parameters:
$I=A/(2\pi)$, where $A$ is the area on the phase portrait of the fast system surrounded by a phase trajectory
passing through the point $(p,q)$. For all $p,q$ on the Hamiltonian level line $\{p,q: \,H(p,q,y,x)=h\}$ in
the domain $G_i$ the function $I$ takes the same value $I_i(h,y,x)$.

Introduce the continuous function $\hat I(p,q,y, x)$ as follows:

\begin{tabular}{ll}
$\hat I(p,q,y, x) = I(p,q,y, x),\quad$    & if $(p,q)\in G_{1,2}(y, x);$\\

$\hat I(p,q,y, x) = I(p,q,y, x)/2,\quad$  & if $(p,q)\in G_3(y, x);$    \\

$\hat I(p,q,y, x) = S(y, x)/(2\pi),\quad$ & if $(p,q)\in l_1 \cup l_2 \cup C.$
\end{tabular}

\begin{figure}[ht]
\hspace*{33mm}\epsfig{file=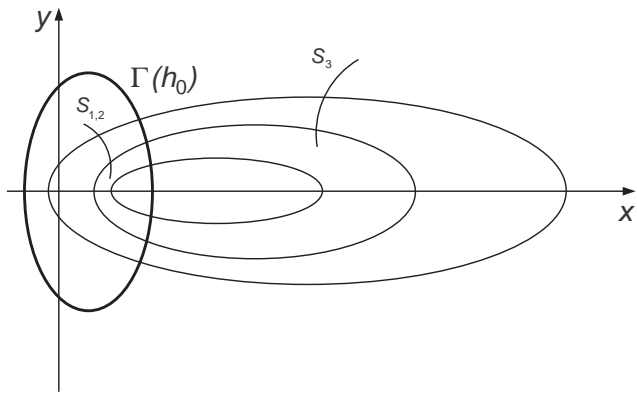,width=90mm} \vspace*{-5mm} \caption{Schematic picture of
motion on the slow plane.}
\label{slowplane}
\end{figure}

The function $\hat I$ is an adiabatic invariant of the perturbed system. Far from the
separatrix of the unperturbed system, its value along a phase trajectory is preserved
with an accuracy of order $\eps$ on time intervals of order $1/\eps$ \cite{Arn63}. When a
phase point crosses a narrow neighborhood of the separatrix, $\hat I$ varies by a value
of order $O(\eps )$ for the majority of initial conditions \cite{N87}. The approximation
$\hat I = \const$ along a phase trajectory is called adiabatic. In this approximation the
dynamics of the slow variables is described by a Hamiltonian system with Hamiltonian
$H_0(\hat I, y,x)$, where $H_0$ is the function $H$ expressed in terms of $\hat I, y,x$.
Consider the dynamics on an energy level $H=h_0$. On the phase plane of the slow
variables $(y,x)$, the separatrix is represented by a curve $\Gm(h_0) = \{y,x:
h_s(y,x)=h_0\}$ (see figure \ref{slowplane}). This curve is called the uncertainty curve
\cite{Wis}. It divides the plane into two domains:  $D_3$, corresponding to $G_3$ on the
fast plane, and the other, $D_{1,2}$, corresponding to $G_1$ and $G_2$. In the adiabatic
approximation, a phase point on the slow plane moves along level lines $I(h_0,y,x) =
\const$. As the area inside the separatrix loops varies, a phase trajectory can cross the
separatrix. Fix $ \nu=1$ or $\nu=2$. Let $B_{\nu}(\hat I_0)$ be the adiabatic trajectory
with $\hat I=\hat I_0$ on the plane of the slow variables for phase points passing
through $G_{\nu}$:
\be
B_{\nu} (\hat I_0) = \{y,x: I_3(h_0,y,x)=2\hat I_0,\, (y,x)\in D_3;\; I_{\nu}(h_0,y,x)=\hat I_0, \, (y,x)\in
D_{1,2}\}
\label{2.2}
\ee

Let us make the following assumption.

B) There exists $\hat I_*$ such that the curve $B_{\nu} (\hat I_*)$ is closed and intersects the
uncertainty curve $\Gm(h_0)$  at two points; at these points $\Th (x,y) \ne 0$. On $B_{\nu}$ there
are no stationary points of the Hamiltonian $H_0$.

If assumption B) is valid both for $\nu = 1$ and $\nu = 2$ (with the same $\hat I_*$), then in the
adiabatic approximation, the projection of the motion with $\hat I = \hat I_0$ onto the slow plane looks
as follows. In the domain corresponding to $G_3$, a phase point moves along the curve $I_3(h,y,x) =
2\hat I_*$. After crossing $\Gm(h_0)$, it follows one of the curves $I_{\nu}(h_0,y,x) = \hat I_*$
[in general, these two curves are different]. Then the phase point crosses $\Gm(h_0)$ again and
continues its motion along $I_3(h_0,y,x) = 2\hat I_*$.

Assumption B) implies that there exists an interval $\Xi$ of values of $\hat I$, such that for all $\hat I_0
\in D$ the curve $B_{\nu}(\hat I_0)$ has the property B).  In the exact system, multiple separatrix crossings
occur in this interval of values of $\hat I$. Computer simulations demonstrate diffusion of the adiabatic
invariant and chaotic dynamics in this region, see, e.g., \cite{NSTV}.

The main result of the paper can be formulated as follows.

{\bf Theorem.} Under additional generality conditions, in the
domain $\{p,q,y,x:\,\hat I \in \Xi$\} on the energy level $H=h_0$
there exist more than $C_1^{-1}/\eps$ stable periodic trajectories
of system (\ref{2.0}). Each of these trajectories is surrounded by
an invariant torus bounding a domain of (3-dimensional) phase
volume greater than $C_2^{-1}\eps$. Inside each of these domains
the variation of $\hat I$ is smaller than $C_3 \eps$.

{\bf Corollary.} The total measure of these  domains is larger than $C_4^{-1}=C_1^{-1}C_2^{-1}$.

Here $C_i$ are positive constants independent of $\eps$. The generality conditions mentioned in the
theorem are formulated below, in Sections \ref{solutions} and \ref{islands}.

The domains bounded by invariant tori are domains (or islands) of perpetual adiabatic invariance. They are
also called stability islands. Projected on the phase plane of the fast variables
$(p,q)$, the motion along the stable periodic trajectories under consideration looks as follows. A phase
point starts moving in region $G_3$, then crosses the separatrix and gets into region $G_{\nu}$,
then crosses the separatrix again, comes back into $G_3$, and returns to its starting point.

The stable periodic trajectories mentioned in the theorem correspond to stable stationary points of the
return map generated by the original system. This map is written in the principal approximation in Section \ref{returnmap}.

{\bf Remark}. The result of the theorem is valid also if the phase portrait of the fast system is of the
kind shown in figure \ref{fastplane} but not for all $x,y$. It is enough to require that the system has
such a phase portrait at the time of separatrix crossing calculated in the adiabatic approximation.

\section{ Adiabatic and improved adiabatic approximations}\label{approximations}

In the fast system, action-angle variables $I, \, \ffi \,\mod
\,2\pi$ are introduced separately in each of the $G_i$ domains by a
canonical transformation of variables.  The corresponding
generating function $W(I,q,y,x)$ contains $y,x$ as parameters [for
brevity, we omit subscripts $i$]. We take this function in the
form

\be W(I,q,y,x) = \int_{q_0(I,y,x)}^q \calP (I,q',y,x) \dd q',
\label{3.0} \ee where $\calP$ is the value of the $p$-variable along
the trajectory with the prescribed value of the action $I$, and
$q_0(I,y,x)$ defines a curve in $G_i$ transversal to the phase
trajectories. In the new variables the Hamiltonian has the form
$H=H_0(I,y,x)$.

Now make a canonical transformation of variables $(p,q,y, x) \mapsto (\bar I, \bar \ffi, \bar y,\bar x)$ with
the generating function $\bar y \eps^{-1} x + W(\bar I,q,\bar y,x)$. The canonically conjugated pairs of
variables are $(\bar I, \bar \ffi)$ and $(\bar y, \eps^{-1} \bar x)$. Formulas for the transformation of
variables are:

\be
\bar \ffi = \pt W/\pt \bar I, \;\; p = \pt W/\pt q, \;\; \bar x = x + \eps \pt W/\pt \bar y, \;\; y
= \bar y + \eps \pt W/\pt x. \label{3.1}
\ee
In the new variables, the Hamiltonian $H$ has the form

\be
H=H_0(\bar I,\bar y,\bar x)+ \eps H_1(\bar I, \bar \ffi, \bar y,\bar x)+ O(\eps^2),
\label{3.2}
\ee
where
\be
H_1 = \frac{\pt H}{\pt y}\frac{\pt W}{\pt x} - \frac{\pt H_0}{\pt x}\frac{\pt W}{\pt \bar y}.
\label{3.3}
\ee
In the adiabatic approximation, the dynamics is described by the Hamiltonian $H_0$. In this approximation
$\bar I=\const$ along phase trajectories.

One can also construct a canonical, close to the identity, transformation of variables $(\bar I, \bar
\ffi, \bar y,\bar x) \mapsto (J,\psi,\hat y,\hat x)$ in order to make the terms of order $\eps$ in
the Hamiltonian independent of the phase (see, for example, \cite{NV05}). In the new variables, the
Hamiltonian takes the form:

\be
\calH = H_0(J,\hat y,\hat x) +\eps \calH_1(J,\hat y,\hat x) +\eps^2 \calH_2(J,\psi,\hat y,\hat
x,\eps), \;\; \calH_1 = \langle H_1 \rangle, \label{3.5}
\ee
where the brackets denote averaging with respect to $\bar \ffi$.

In the improved adiabatic approximation, the dynamics is described by the Hamiltonian $H_0(J,Y,X) +
\eps \calH_1(J,Y,X)$. In this approximation $J$ is an integral of motion. With an accuracy of
order $\eps^2$, the following formula for $J$ is valid (see \cite{N87}):
\bea
J = J(p,q,y,x) = I + \eps u(p,q,y,x), \label{3.9} \\
\nonumber\\
u = \frac{1}{4\pi}\left[\int_0^T\left(\frac{\pt E}{\pt y}\int_0^t \frac{\pt E}{\pt x} \dd\sigma \right)\dd t
- \int_0^T\left(\frac{\pt E}{\pt x}\int_0^t \frac{\pt E}{\pt y} \dd\sigma \right)\dd t \right] \nonumber
\\+\frac{1}{2\pi} \int_0^T \left(\frac{T}{2} - \sigma \right) \left\{E,h_s\right\} \dd \sigma.
\label{3.9a}
\eea
The integrals here are calculated along a phase trajectory of the fast system passing through the point
$(p,q)$; $\sigma$ is the time of motion along this trajectory starting from this point, $T$ is the period of
motion. The function $J$ is the improved adiabatic invariant. In the complete system far from separatrices its
value along a phase trajectory is constant with an accuracy of order $\eps^2$ on time intervals of order
$\eps^{-1}$. In what follows, it is convenient to use the function $\hat J(p,q,y,x)$ equal to $J$ in the domains
$G_1,G_2$ and equal to $J/2$ in the domain $G_3$.

\section{ Description of the separatrix crossing }\label{crossing}

On the phase plane $\mathbb{P}^{sl}$ of the slow variables $(y,x)$, the separatrix is represented by a curve
$\Gm(h_0)$ (see figure \ref{slowplane}). Let $B_{\nu}(\hat I) \in \mathbb{P}^{sl}$ be a curve (see
(\ref{2.2})) corresponding to $\hat I$ such that $\hat I \in \Xi$. According to B), $B_{\nu}(\hat I)$ is
divided by $\Gm(h_0)$ into two parts. Consider two segments $\calS_3$ and $\calS_{\nu}$ on $\mathbb{P}^{sl}$
separated by $\Gm(h_0)$, such that each of the segments crosses transversely all the curves $B_{\nu}(\hat
I)$ with $\hat I \in \Xi$ at a distance of order 1 from $\Gm(h_0)$. Points on the fast plane corresponding to
$\calS_3$ belong to $G_3$, and the points corresponding to $\calS_{\nu}$ belong to $G_{\nu}$. Fix an interval
$\Xi_0 \in \Xi$, assuming that the  endpoints of $\Xi_0$ and $\Xi$ are different.

Let the motion start at $t=0$ in a point $M^{(0)}(p^{(0)},q^{(0)},y^{(0)}, x^{(0)})$,
such that $H(p^{(0)},q^{(0)},y^{(0)}, x^{(0)})=h_0$ and on the plane $\mathbb{P}^{sl}$
$(\hat y^{(0)},\hat x^{(0)}) \in \calS_3$. Hence, $(p^{(0)},q^{(0)})\in G_3(y^{(0)},
x^{(0)})$. Let $\hat I = \hat I^{(0)}, \; \hat J=\hat J^{(0)}, \; \psi = \psi^{(0)}$ at
this point; $\hat I^{(0)}\in \Xi_0$. In the adiabatic approximation the trajectory of the
slow variables is $B(\hat I^{(0)})$, the motion is described by a Hamiltonian system with
Hamiltonian $H_0(\hat I^{(0)}, y,x)$. Let $(x_-, y_-)$ and $(x_+, y_+)$ denote the points
of intersection of $B(\hat I^{(0)})$ and $\Gm(h_0)$; the sign ``$-$'' corresponds to
passage from $G_3$ to $G_{\nu}$ and the sign ``+'' corresponds to passage from $G_{\nu}$
to $G_3$. Below $\tau_-$ and $\tau_+$ denote the slow time moments of separatrix crossing
in this approximation. In the exact system, in the process of evolution the phase point
on $\mathbb{P}^{sl}$ approaches the curve $\Gm(h_0)$, and accordingly, on the plane
$(p,q)$ it approaches the separatrix. We are interested in the dynamics of phase points
that are being captured into $G_{\nu}$ after the separatrix crossing. For such a phase
point, the first crossing of the $Cq$-axis in $G_{\nu}$ (see figure \ref{fastplane})
occurs near $C$ at time $\tau = \tau_- + O(\eps)$. [We assume for simplicity that
coordinate axes in figure \ref{fastplane} are the principal axes of the saddle point.]
Denote the value of $E$  at the point of this crossing as $h^{(0)}$. We introduce
$\eta^{(0)} = 1- |h^{(0)}/\eps\Th_-|$\,. Here $\Th_- \equiv \Th (y_-, x_-) > 0$.

After the separatrix crossing the phase point in the plane of slow variables $(\hat y,\hat x)$ moves towards
$\calS_{\nu}$. When it crosses $\calS_{\nu}$, the projection of the phase point onto the fast plane is deep
inside region the $G_{\nu}$. Denote the value of $\hat J$ at this time moment as $\hat J^{(1)}$.

Then the phase point starts again approaching the separatrix.  At $\tau = \tau_+ +
O(\eps)$ the phase point crosses $Cq$-axis in $G_{\nu}$ near the point $C$ for the last
time before entering $G_3$. Denote the value of $E$  at the point of this crossing as
$h^{(1)}$. We introduce $\eta^{(1)} =|h^{(1)}/\eps\Th_+|$\,. Here $\Th_+ \equiv \Th(y_+,
x_+) < 0$.

After crossing $\Gm(h_0)$, the projection of the phase point onto the plane $(\hat y,\hat
x)$ crosses again the segment $\calS_3$. Let $\hat J^{(2)}, \;  \psi^{(2)}$ denote the
values of $\hat J, \; \psi$ at this time moment. Then the phase point approaches the
separatrix again, crosses it and gets captured into $G_l,\, l=1$ or $l=2$. Let
$E=h^{(2)}$ at the first crossing of the $Cq$-axis in $G_l$ near $C$. This crossing
occurs at time $\tau = \tau_- + T_0 + O(\eps)$, where $T_0$ is the slow time period of
motion along the trajectory $\hat I=\hat I^{(0)}$ in the adiabatic approximation. We
introduce $\eta^{(2)} = 1- |h^{(2)}/\eps\Th_-|$\,. We are interested in the dynamics of
the phase points for which $l=\nu$.

\section{ The return map}\label{returnmap}

In the energy level $H=h_0$ the segment $\calS_3$ is represented by a piece of a two-dimensional surface
$\{p,q,y,x: H(p,q,y,x) = h_0, (\hat y,\hat x) \in \calS_3 \}$. This piece can be parametrized by variables
$\hat J, \psi$. The corresponding Poincar\'e return map $M\, :\,(\hat J^{(0)},  \psi^{(0)}) \to (\hat
J^{(2)},  \psi^{(2)})$ produced by trajectories that pass through $G_{\nu}$ is symplectic. Its stable
stationary points correspond to stable periodic orbits of the original problem, of period approximately equal
to $T_0/\eps$. It is convenient to study these stationary points using a different set of variables, namely
to consider the map $\hat M\, :\, (\hat J^{(0)},\eta^{(0)},\nu) \to (\hat J^{(2)},\eta^{(2)},l)$.

The map $\hat M$ is the composition of two maps: $\hat M=M^{(2)}\circ M^{(1)}$,
$$
M^{(1)}:(\hat{J}^{(0)},\eta ^{(0)},\nu)\to (\hat{J}^{(1)},\eta ^{(1)},\nu)\,,\quad
M^{(2)}:(\hat{J}^{(1)},\eta ^{(1)},\nu)\to (\hat{J} ^{(2)},\eta ^{(2)},l)\,.
$$

The results of \cite{N87,NV05} give the following formulas for $M^{(k)}$. Suppose that
\be
\eta ^{(0)},\eta ^{(1)},\eta ^{(2)}\in [c_1^{-1},1-c_1^{-1}]\,.
\label{5.1}
\ee
Then
\smallskip
\bea
\hat{J}^{(1)}  &=&   \hat{J}^{(0)}-(2\pi )^{-1}\eps a_{-} \Theta_{-}\ln (2\sin \pi \eta
^{(0)})+(2\pi)^{-1}\eps (d_{\nu,-}-d_{3,-}/2)+ \nonumber \\
& & \quad (2\pi)^{-1}\eps\Th_-(\eta^{(0)} - 1/2)(b_{\nu,-} -b_{3,-}/2)+O(\eps ^{3/2}\ln \eps )\,,
\label{5.2f}\\
\eta ^{(1)}  &=&  \{\eta ^{(0)}+\eps ^{-1}\Phi _1^{(\nu)}(\hat{J}%
^{(1)},\eps )+O(\eps ^{1/3}\ln ^{-1/3}\eps )\}\,,
\label{5.2}
\eea

\bea
\hat{J}^{(2)}  &=&  \hat{J}^{(1)}-(2\pi )^{-1}\eps a_{+}\Theta_+ \ln (2\sin \pi \eta ^{(1)}) -\eps
(2\pi )^{-1}(d_{\nu,+}-d_{3,+}/2) + \nonumber \\
& &\quad (2\pi )^{-1}\eps\Th_+ (1/2 - \eta^{(1)} )(b_{\nu,+} - b_{3,+}/2) + O(\eps ^{3/2}\ln \eps)\,,
\label{5.3f} \\
\eta ^{(2)}  &=&  \{\eta ^{(1)}+\eps ^{-1}\Phi _2(\hat{J} ^{(2)},\eps )+O(\eps ^{1/3}\ln ^{-1/3}\eps
)\}\,,\label{5.3}
\eea

\bea
l=\nu,\quad  & \mbox{if} & \quad 0<\{(\eps ^{-1}\Phi _2(\hat J^{(2)}, \eps)+\eta ^{(1)})/2\}<1/2,
\label{lnu1}\\
l\ne \nu,\quad & \mbox{if} & \quad 1/2<\{(\eps ^{-1}\Phi _2 (\hat J^{(2)}, \eps)  +\eta ^{(1)})/2\}<1.
\label{lnu2}
\eea
\smallskip

Here $\{\,\cdot \, \}$ denotes the fractional part, \bea a_{\pm
}=a(y_{\pm},x_{\pm}),\quad b_{j,\pm} = b_j (y_{\pm},x_{\pm}),
\quad
d_{j,\pm}=d_j(y_{\pm},x_{\pm}), \nonumber\\
j=1,\, 2,\, 3, \quad b_3=b_1+b_2, \quad d_3=d_1+d_2. \nonumber
\eea
For $a$ one has $a=1/\sqrt{-g } $, where $g $ is the Hessian of~$E$ at the point $C$. The values $d_{1,2}$ are
defined as the main terms of the expansions (cf. \cite{N87})
$$
2\pi u(p,q,y,x)=d_j+O(\sqrt{|E|} \ln |E| ),
$$
where $u$ is given by (\ref{3.9a}), and the point $(p,q)$ is on the $Cq$-axis (figure \ref{fastplane}) in
$G_j$ near $C$. Values $b_j$ come from the expansion of the period of motion in $G_j, \; j=1,2$ near the
separatrix:
$$
T_j = -a\; \ln|E| + b_j + O(E\ln|E|).
$$
In (\ref{5.2}), (\ref{5.3}),
    \smallskip
\bea
\Phi _1^{(\nu)}(\hat J,\eps )&=&{\frac{{1}}{{2\pi }}}\int_{\tau _{-}^{(\nu)}}^{\tau _{+}^{(\nu)}}(\omega
_0^{(\nu)}(\hat J,Y_{\nu} (\tau),X_{\nu} (\tau))+\eps \omega _1^{(\nu)}(\hat J,Y_{\nu} (\tau),X_{\nu}
(\tau)))\,\dd\tau \,,
\label{5.3a}\\
 \Phi _2(\hat J,\eps )&=&{\frac{{1}}{{\pi }}}\int_{\tau _{+}^{(3)}}^{\tau  _{-}^{(3)}}(\omega _0^{(3)}(2\hat J,Y_3
 (\tau),X_3 (\tau))+  \eps \omega _1^{(3)}(2\hat J,Y_3 (\tau),X_3 (\tau)))\,\dd\tau .
 \label{5.3b}
 \eea
\smallskip
Here $\omega _0^{(j)}$  and  $\omega _1^{(j)}$ are defined as $\omega _0^{(j)}=\partial H_0/\partial J$,
$\omega _1^{(j)}=\partial \calH_1/\partial J$, with $H_0,\calH_1$ calculated in the region $G_j$, and $(Y_j,X_j)$
is a solution of the Hamiltonian system with Hamiltonian $\hat H = H_0(J,y,x) + \eps H_1 (J,y,x) $ on the energy
level $\hat H=h_0$. The values $\tau_{\pm}^{(j)}$ are the slow time moments when the phase point corresponding to
this solution arrives to the separatrix. At these moments $ S_j (Y,X) =2\pi J$. [One can take any of such
solutions, they differ by a time shift which does not change the value of the integrals in (\ref{5.3a}),
(\ref{5.3b})]. Conditions (\ref{5.1}) are understood inductively: if $\eta ^{(0)}\in [c_1^{-1},1-c_1^{-1}]$,
and $\eta ^{(1)}$, found according to (\ref{5.2f}), (\ref{5.2}) lies in segment $[c_1^{-1},1-c_1^{-1}]$, then
(\ref{5.2}) is valid, and analogously for $\eta ^{(2)}$ and expressions (\ref{5.3f}), (\ref{5.3}).

Equations (\ref{5.2f}), (\ref{5.3f}) follow directly from the formula for the jump of the adiabatic invariant at a
separatrix \cite{N87} for the case when the areas of domains $G_1$ and $G_2$ are equal. Equation (\ref{5.2})
follows directly from the formula for phase change between separatrix crossings \cite{NV05}. Equations
(\ref{5.3})-(\ref{lnu2}) follow from the result of \cite{NV05} in the following way.

Let $t_+, h_+$ be the values of $t$ and $E$ when the phase point
crosses the axis $Cp$ near $C$ for the first time after exit from
$G_{\nu}$. At time $t_+$ the phase point is in $G_3$ on the positive
or negative part of the axis $Cp$. Let $h_-$ be the value of $E$ when
the phase point crosses the same part of $Cp$ for the last time
before exit from $G_3$. Introduce notations: $\bar\zeta =
\frac{h_+}{2|\Th_+|}, \; \zeta = \frac{h_-}{2|\Th_-|}$. According
to \cite{NV05},
\be
\bar\zeta + \zeta = \frac{1}{2\eps} \Phi
_2(\hat{J} ^{(1)},\eps )+O(\eps ^{1/3}\ln ^{-1/3}\eps )\; \mod
\;1,
\label{5.5}
\ee
provided $\bar\zeta,\zeta \in [c_2^{-1},1-c_2^{-1}]$. The values $h^{(1)}$ (see Section
\ref{crossing}) and $h_+$ are related as $h_+ = h^{(1)} + \eps |\Th_+| + O(\eps^{3/2})$
\cite{N87}. Therefore, $\bar\zeta = \frac12 (1-\eta^{(1)}) + O(\eps^{1/2})$. The values
$h^{(2)}$ (see Section \ref{crossing}) and $ h_-$ are related as $ h_- = h^{(2)} + 2\eps
\Th_- + O(\eps^{3/2})$ if  $l=\nu$, and as $ h_- = h^{(2)} + \eps \Th_- + O(\eps^{3/2})$
if  $l\ne\nu$. We have $l=\nu$ provided that $1/2 < \zeta < 1$. Consider this case. We
obtain that $\zeta = \frac12(1 + \eta^{(2)}) + O(\eps^{1/2})$. From (\ref{5.5}) we find
\be
\frac12 (1-\eta^{(1)}) + \frac12 (1 + \eta^{(2)}) =
\frac{1}{2\eps} \Phi _2(\hat{J} ^{(2)},\eps )+O(\eps ^{1/3}\ln ^{-1/3}\eps ) \; \mod \;1.
\nonumber
\ee
Therefore,
\be
\eta^{(2)} = \eta^{(1)} + \frac{1}{\eps} \Phi _2(\hat{J}
^{(2)},\eps )+O(\eps ^{1/3}\ln ^{-1/3}\eps ) \; \mod \;2.
\label{5.6}
\ee
Then condition $1/2 <\zeta<1$ implies equation
(\ref{lnu1}).

Now consider the case $l \ne \nu$. In this case $0<\zeta<1/2$. We get that $\zeta$ and
$\eta^{(2)}$ are related as $\zeta=\frac12 \eta^{(2)} + O(\eps^{1/2})$. From (\ref{5.5})
we find
\be
\frac12 (1-\eta^{(1)}) + \frac12 \eta^{(2)} = \frac{1}{2\eps} \Phi _2(\hat{J} ^{(2)},\eps )+O(\eps
^{1/3}\ln ^{-1/3}\eps ) \; \mod \;1. \nonumber
\ee
Therefore,
\be
\eta^{(2)} = \eta^{(1)} + \frac{1}{\eps} \Phi _2(\hat{J} ^{(2)},\eps ) -1 +O(\eps ^{1/3}\ln
^{-1/3}\eps ) \; \mod \;2.
\label{5.7}
\ee
Combining (\ref{5.6}) and (\ref{5.7}) we obtain equation (\ref{5.3}). Condition $0<\zeta<1/2$ implies
(\ref{lnu2}).

\section{Linearly stable solutions}\label{solutions}

\subsection{Equations for stable stationary points}

Linearly stable periodic solutions of the original problem of period approximately equal to $T_0/\eps$
correspond to linearly stable fixed points of the map $\hat M$. These points are defined by the following set
of conditions:
\bea
\hat J^{(2)} = \hat J^{(0)}, \quad  \eta^{(2)} = \eta^{(0)},\quad l=\nu,
\label{6.1} \\
|\,\mbox{tr}\, \hat M'\,|<2,
\label{6.2}
\eea
where $\hat M'$ is the linearization of $\hat M$ at the fixed point.

Put $\Phi_1 \equiv  \Phi_1^{(\nu)}$. Introduce also the following notations:
\bea
\hat J^{(s)} = \eps\xi^{(s)}\,,\quad  s = 0,1,2;
\label{6.3} \\
\eps^{-1} \frac {d\Phi_k(\eps\xi,\eps)}{d\xi} = \gamma_k(\eps\xi,\eps)\,,\qquad
\gamma_k(\eps\xi,0)=\gamma_k^0(\eps \xi)\,, \quad k=1,2;
\nonumber \\
a_{\pm }\Theta_\pm = \mp 2\pi\alpha_\pm\,,\quad d_{\nu,\pm}-d_{3,\pm}/2=- 2\pi d_\pm \,,\quad
\Th_{\pm}(b_{\nu,\pm}-b_{3,\pm}/2) = 2\pi b_{\pm}. \nonumber
\eea
Suppose that $(\xi^{(0)},\eta^{(0)})=(\xi ,\eta )$ correspond to a stationary point of the map $\hat M$. We
rewrite equations (\ref{6.1}) using formulas (\ref{5.1})-(\ref{5.3}), notations (\ref{6.3}), and neglecting
terms $O(\eps^{1/3}\ln^{-1/3}\eps )$:
\bea
\alpha_- \ln (2\sin\pi\eta)\!+\!d_- \!-\!  (\eta \!-\! 1/2)b_- = \alpha_+ \ln(2\sin\pi\eta^{(1)})\!+\!d_+
\!+\!(1/2 \!-\! \eta^{(1)})b_+,
\label{6.4}\\
\eta^{(1)}-\eta +\gamma_1^0\alpha_-\ln (2\sin\pi \eta)+\gamma_1^0d_- - \gamma_1^0 (\eta - \frac12)b_-
=\eps^{-1}\Phi_1(\eps \xi,\eps )\, \mod\, 1  ,
\label{6.5}\\
\eta -\eta^{(1)}=\eps^{-1} \Phi_2(\eps \xi ,\eps )\, \mod\, 1 ,
\label{6.6}\\
\eta^{(1)}+\eps^{-1} \Phi_2(\eps \xi ,\eps ) \in (0,1)\, \mod\, 2 ,
\label{6.7}\\
\eta,\eta^{(1)} \in (c_1^{-1}, 1 - c_1^{-1}). \label{6.8}
\eea

The stability condition $|\,\mbox{tr}\, \hat M'\,|<2$ can be written using (\ref{5.2f})-(\ref{5.3}) [we again
neglect terms $O(\eps^{1/3}\ln^{-1/3}\eps )$]:
\bea
-4<Q<0 , \label{6.9} \\
Q = (u_2 - u_1)(\gamma_2^0 + \gamma_1^0) - u_2u_1\gamma_2^0\gamma_1^0 ,
\label{6.10}\\
u_1 = \alpha_- \pi \cot \pi \eta -  b_- ,  \quad  u_2 = \alpha_+ \pi \cot \pi \eta^{(1)} - b_+.
\nonumber
\eea
The exact set of conditions for a stable periodic solution differs from (\ref{6.4})-(\ref{6.8}) by
terms $O(\eps^{1/3}\ln^{-1/3}\eps )$ that can be differentiated in $\eps \xi \,, \eta \,,
\eta^{(1)}$ without changing their smallness order. Below, we find non-degenerate solutions of
system (\ref{6.4})-(\ref{6.8}). According to the implicit function theorem, for small enough $\eps$,
each of these solutions corresponds to a solution of the exact system.

\subsection{Existence of stable stationary points}\label{existence}

We assume that the following generality conditions are valid.

{\bf H1}. In $\Xi_0$ there exists an interval of values of $I$ such that
\be
\gamma_1^0(I)\cdot \gamma_2^0(I)\ne 0\,, \qquad  {\it \frac \dd{\dd
I}}(\gamma_1^0(I)/\gamma_2^0(I))\ne 0\,. \label{6.11a}
\ee
In what follows, we consider only those values of $\eps \xi$ that belong to the interval defined by {\bf
H1}.

Equation (\ref{6.4}) has the form
$$
F(\eta; \al_-, b_-, d_-) = F(\eta^{(1)}; \al_+, b_+, d_+)
$$
with $F(\eta; \al, b, d)= \al \ln (2\sin \pi \eta) +d-b(\eta-1/2), \, \eta\in (0,1),\, \al>0$. We have $F \to
-\infty$ as $\eta \to 0^+$ or $\eta \to 1^-$. Moreover, $F$ has a unique maximum $F_{max} = m(\al, b, d)$ on
the interval $0<\eta<1$. We put $m_{\pm} = m(\al_{\pm}, b_{\pm}, d_{\pm})$.

\begin{figure}[ht]
\begin{tabular}{ccc}
\hspace*{10mm}\epsfig{file=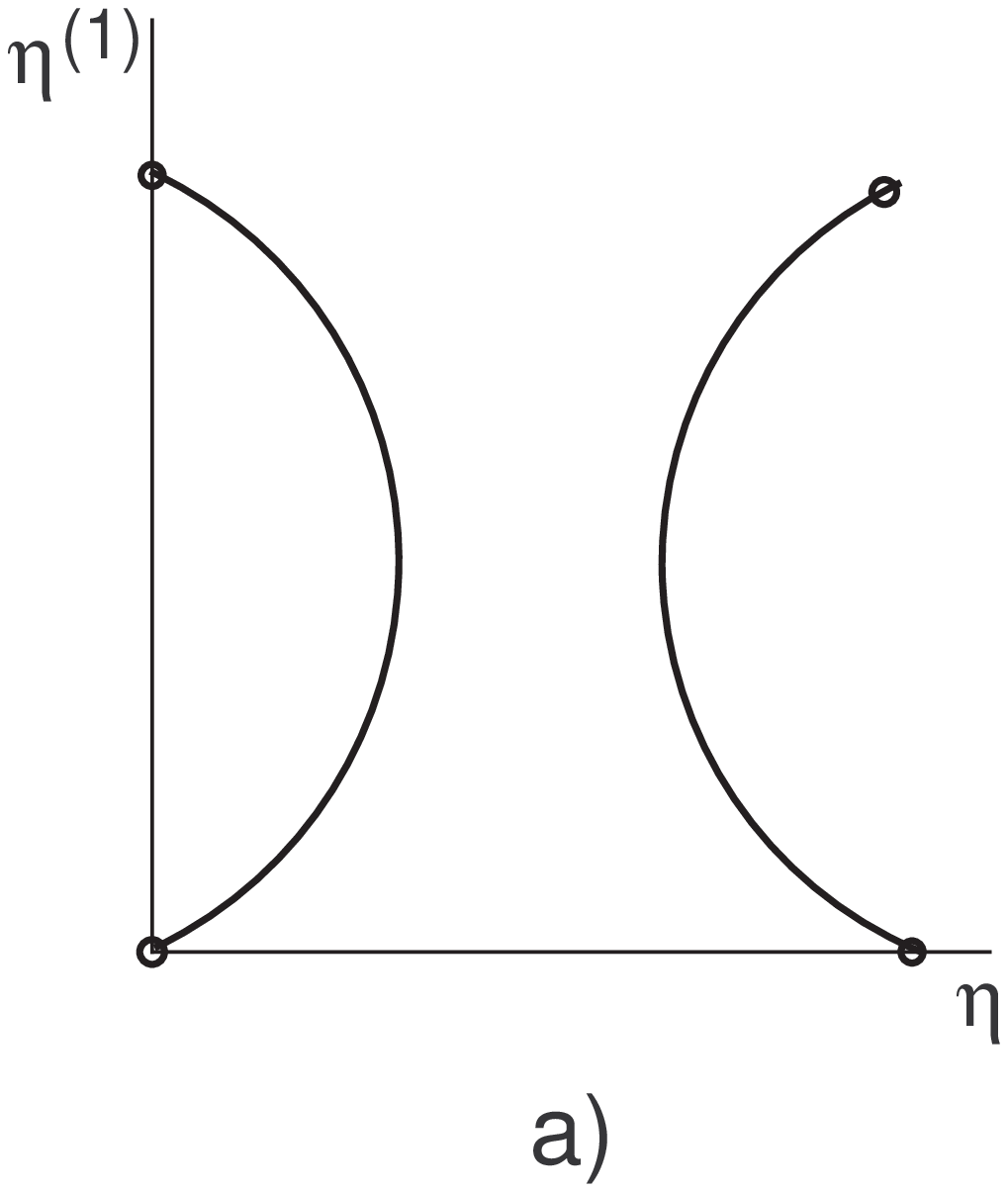,width=40mm}& \hspace*{0mm}\epsfig{file=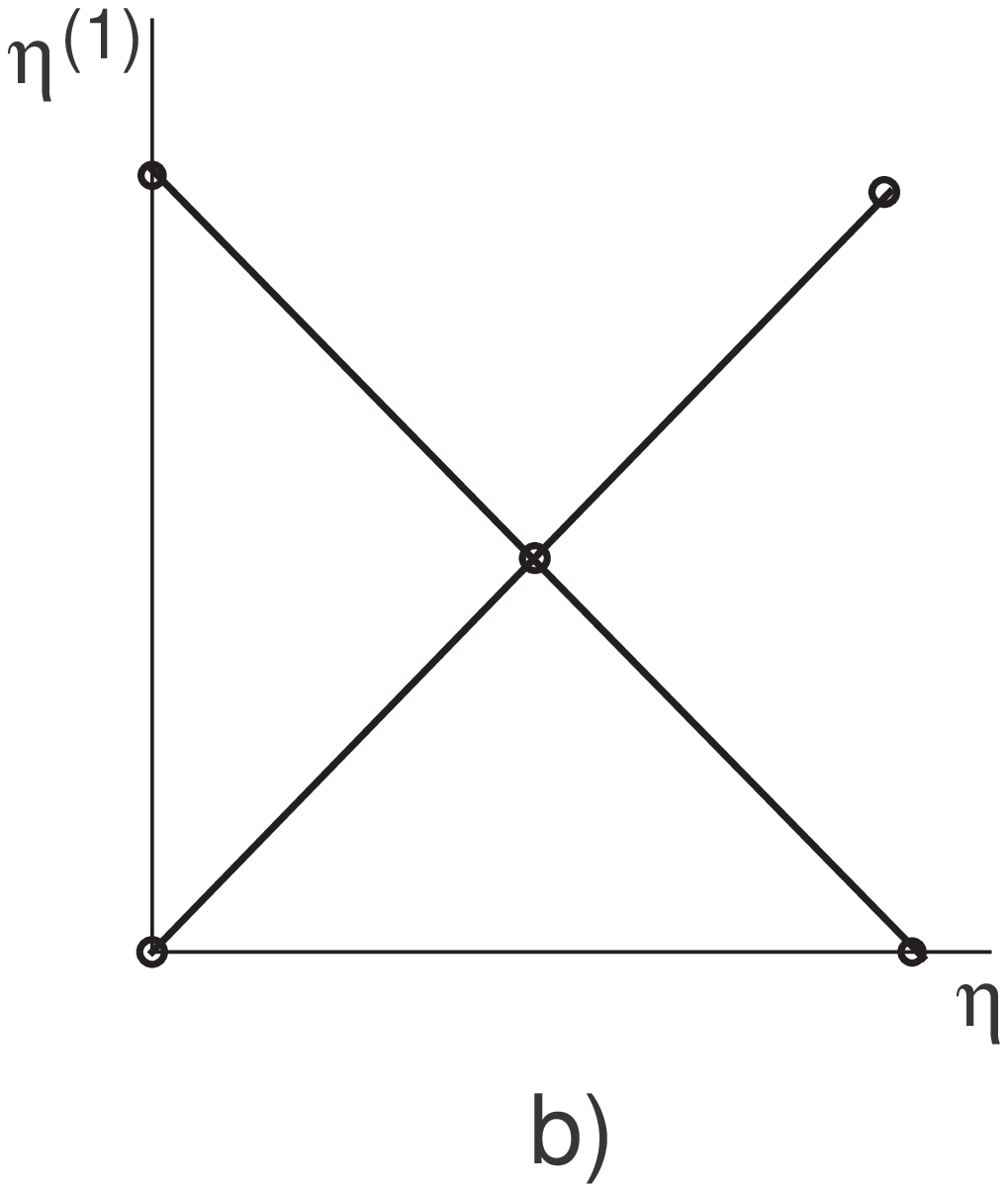,width=40mm}&
\hspace*{-0mm}\epsfig{file=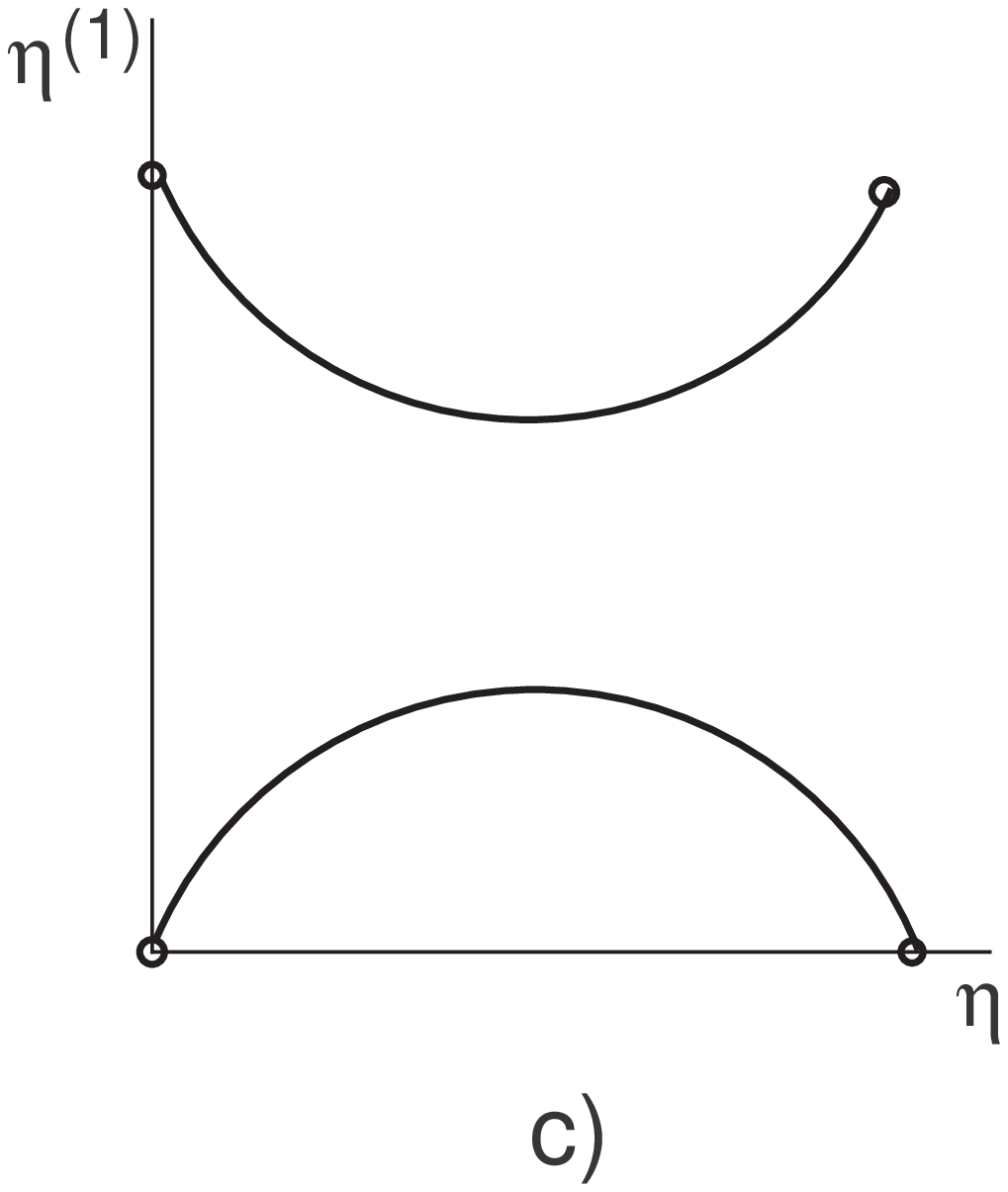,width=40mm}
\end{tabular}
\caption{Set of points $(\eta ,\eta^{(1)})$ satisfying (\ref{6.4}). a) $m_- >m_+$; b) $m_- =m_+$; c)
$m_- <m_+$.}
\label{points}
\end{figure}

The set of points $(\eta ,\eta^{(1)})\in (0,1)\times (0,1)$ satisfying (\ref{6.4}) is shown in
figure \ref{points}. Consider for definiteness the case $m_- >m_+,\quad \gamma_1^0\gamma_2^0<0$. The
other cases can be studied analogously.

Let in figure \ref{points} $A$ and $B$ have coordinates $(0, 0 )$
and $(0,1)$ respectively. Consider condition (\ref{6.9}) on the
curve $AB$. Near the point $A$ both $u_1$ and $u_2$ are large
positive, and therefore $Q$ is large positive. Near the point $B$ both
$u_1$ and $-u_2$ are large positive, and therefore $-Q$ is large
positive. Therefore, on the segment $AB$ there is an interval
$L=L(\eps \xi )$ (or several intervals) where condition
(\ref{6.9}) is valid. If $c_1$ is sufficiently large, the relations
(\ref{6.8}) are also valid on $L(\eps \xi )$.

Now we look for solutions of the system
\bea
 (\eta ,\eta^{(1)})\in  L(\eps \xi ),    \label{6.11}\\
 \eta^{(1)}-\eta + \gamma_1^0\alpha_ - \ln ( 2\sin \pi \eta  ) + \gamma_1^0d_- -\gm_1^0 b_-(\eta - 1/2) =
\eps^{-1}\Phi_1(\eps \xi , \eps ) \quad \mod \, 1, \label{6.12}\\
 \eta - \eta^{(1)} = \eps^{-1}\Phi_2(\eps \xi ,\eps ) \quad \mod \, 2. \label{6.13}
\eea
Relations (\ref{6.11})-(\ref{6.13}) are equivalent to (\ref{6.4})-(\ref{6.10}). Relation (\ref{6.7}) holds
because, according to (\ref{6.13}), $\eta^{(1)} + \eps^{-1}\Phi_2(\eps \xi ,\eps )= \eta \quad \mod \, 2 $
and $\eta \in (0,1)$\,. Validity of the other relations (\ref{6.4})-(\ref{6.10}) follows from the definition
of $L(\eps \xi )$ and equations (\ref{6.12}), (\ref{6.13}).

System (\ref{6.11})-(\ref{6.13}) can be interpreted geometrically as follows. Consider the curve
$\La(\eps\xi)\subset \mathbb{T}^2 = \{(s_1\, \mod \, 1,s_2\, \mod \, 2)\}$, which is the image of $L(\eps \xi
)$ under the mapping
$$
( \eta ,\eta^{(1)})\mapsto (s_1,s_2)=(\eta^{(1)}-\eta +\gamma _1^0\alpha _{-}\ln (2\sin \pi \eta)  +
\gamma^0_1d_- -\gm_1^0 (\eta - 1/2)b_-, \quad \eta -\eta^{(1)})
$$
defined by the left hand sides of (\ref{6.12}), (\ref{6.13}). Consider also the point
$$
\eps ^{-1}\Phi (\eps \xi
,\eps ) = (\eps ^{-1}\Phi _1(\eps \xi ,\eps ) \, \mod \, 1, \quad \eps^{-1} \Phi_2(\eps \xi , \eps )\, \mod
\, 2)
$$
on $\mathbb{T}^2$. As $\xi$ varies, the curve $\La(\eps\xi)$ is moving slowly on the torus. At the same time,
the point $\eps ^{-1}\Phi$ is moving fast, with its velocity vector $\gamma (\eps \xi , \eps )=(\gamma
_1(\eps \xi ,\eps )\,, \gamma _2(\eps \xi ,\eps ))$ varying slowly. Solutions of system
(\ref{6.11})-(\ref{6.13}) correspond to values of ``time'' $\xi$ when point $\eps ^{-1}\Phi(\eps \xi ,\eps )$
crosses the curve $\La(\eps\xi)$.

It follows from assumption {\bf H1}, that there exists $I_c$ such that the ratio
$\gamma_1^0(I_c)/\gamma_2^0(I_c)$ is irrational. Then in the interval $\xi \in (I_c-c_3^{-1}/\eps ,\,
I_c+c_3^{-1}/\eps )$ on every segment of length $c_4$ of ``time'' $\xi$ the point $\eps ^{-1}\Phi(\eps \xi
,\eps )$ crosses transversely the moving curve $\La(\eps\xi)$ at least once. Indeed, the point's trajectory
is close to a dense winding of the torus with irrational angular coefficient
$\gamma_1^0(I_c)/\gamma_2^0(I_c)$, and the curve $\La(I_c)$ is transversal to the vector
$\gm^0(I_c)=(\gamma_1^0(I_c), \gamma_2^0(I_c))$ at least at one point. [The latter is because $\La(I_c)$ is
not straight.]

Therefore, there exist more than $c_3^{-1}c_4^{-1}/\eps -1 >C_1^{-1}/\eps $ intersection points. Each
intersection point corresponds to a linearly stable stationary point of the return map. In Section
\ref{distrib} we obtain an asymptotic expression for the number of intersection points.

\section{Stability islands}\label{islands}

Let $\xi_r,\, \eta_r,\, \eta_r^{(1)}$ be one of the solutions of system
(\ref{6.11})-(\ref{6.13}), found in the previous section. One can rewrite the map $\hat
M$ in variables $\tilde\xi  = \xi -\xi_r,\, \eta $, expand its right hand sides in $\eps$
supposing $\tilde\xi \sim 1$, and neglect terms $O( \eps^{1/3}\ln ^{-1/3}\eps)$. Thus,
one obtains a map that does not contain $\eps$ and depends on the parameter $I_r = \eps
\xi_r$. This map has a stationary linearly stable point $\tilde\xi =0,\, \eta=\eta_r$.
According to Kolmogorov-Arnold-Moser (KAM) theory, under certain conditions [absence of
resonances up to the 4-th order, non-zero coefficient in the normal form] the stationary
point of the complete map is surrounded with a Cantor family of invariant curves (see,
e.g. \cite{AKN}). These latter conditions of KAM theory in the case under consideration
can be written as one condition, namely $f(I_r,\, \eta_r) \ne 0$, where $f$ is a
continuous function. We assume that the following condition is valid.

{\bf H2}. There exist $I_*,\, \eta_*,\, \eta_*^{(1)}, \quad ( \eta_*,\, \eta_*^{(1)} ) \in L(I_*)$ such that
$f(I_*,\, \eta_*) \ne 0$\,.

Then for $I$ close to $I_*,\, |I-I_*| < c_5^{-1}$, there is a segment of the curve $L(I)$ where
$|f(I,\,\eta,\, \eta^{(1)})| > c_6^{-1}$. Without loss of generality one can assume that this segment
coincides with $L(I)$ [otherwise, we can shorten $L(I)$]. According to KAM theory, for $|I_r-I_*| < c_5^{-1}$
there is an invariant curve at distance bigger than $c_6^{-1}$ from the point $\tilde\xi=0,\, \eta=\eta_r$.
We call the domain surrounded by this curve a stability island of the Poincar\'e map. The area of this island
is larger than $c_7^{-1}$. Returning from $\tilde\xi ,\, \eta $ to the original variables $\xi,\, \eta $, we
find that the area of the stability island is estimated from below as $c_7^{-1}\eps$.

The boundary of the stability island (an invariant curve) is a section of an invariant torus on the energy
level of the original system. The volume of the domain surrounded by this torus is larger than $C_2^{-1}
\eps$. The variation of $\hat I$ inside this domain is smaller than $C_3 \eps$.

\section{Distribution of periodic orbits}\label{distrib}

Let $N(I_c, \Delta, \eps)$ be the number of stable periodic orbits under consideration with $I^{(0)} \in
[I_c, \, I_c +\Dt] \subset \Xi_0$. According to Section \ref{existence} this number is equal to the number of
``time'' points $\xi \in [I_c/\eps, \, (I_c +\Dt)/\eps]$ when the point $\eps^{-1} \Phi(\eps \xi, \eps)$
crosses the curve $\La(\eps\xi)$ on the torus $\mathbb{T}^2 = \{(s_1\, \mod \, 1,s_2\, \mod \, 2)\}$. For
values of $I_c$ such that $\gm_1^0(I_c)/\gm_2^0(I_c)$ is irrational one can calculate $\rho (I_c) =
\lim_{\Dt \to 0}\lim_{\eps \to 0} (\eps N(h_0,I_c)/\Dt)$. We give the result omitting the standard argument
based on ergodicity of winding of a torus with frequencies $(\gm_1^0(I_c), \, \gm_2^0(I_c))$ (cf.
\cite{Tal}). Let $l$ be a natural parameter on $\La(I_c)$, and $n(l,I_c)$ be a unit normal vector on
$\La(I_c)$. Introduce the notation
\be
\chi (I_c) = \int_{\La(I_c)} |(\gm^0 (I_c), \, n(l,I_c))|\, \dd l,
\ee
where $(\gm^0, \, n)$ is the Euclidean scalar product of $\gm^0$ and $n$. Hence, $\chi$ is an ``absolute''
phase flux of the vector field $\gm^0$ across the curve $\La$: for the phase flux across any element of $\La$ we
take its absolute value. The value of $\chi (I_c)$ equals the product of the length of vector $\gm^0$ and the
full length of the projection of the curve $\La(I_c)$ onto the normal to vector $\gm^0$. Then
\be
\rho (I_c) = \frac12 \chi (I_c)
\ee
(the coefficient $1/2$ in this formula appears because the surface area of the torus equals 2).

Let $\bar N(\Xi_0, \eps)$ be the number of stable periodic orbits under consideration with $I^{(0)}
\in \Xi_0$. If $\gm_1^0 \dd \gm_1^0/\dd I - \gm_2^0 \dd \gm_2^0/\dd I$ vanishes only on a set of
measure 0 on  $\Xi_0$, then
\be
\lim_{\eps \to 0} \eps \bar N(D_0,\eps) = \frac12 \int_{\Xi_0} \chi (I) \dd I.
\ee

\section*{Acknowledgements}

The work was supported in part by grants RFBR 06-01-00117, 05-01-01119 and NSh-1312.2006.1.
C.S. has been supported by grants BFM2003-09504-C02-01, MTM2006-05849/Consolider
(Spain) and CIRIT 2005 SGR-1028 (Catalonia).

\end{document}